%
%
\input amstex.tex
\documentstyle{amsppt}
\magnification=1200
\baselineskip=13pt
\hsize=6.5truein
\vsize=8.9truein
\parindent=20pt
\newcount\sectionnumber
\newcount\equationnumber
\newcount\thnumber
\newcount\countrefno

\def\ifundefined#1{\expandafter\ifx\csname#1\endcsname\relax}
\def\assignnumber#1#2{%
        \ifundefined{#1}\relax\else\message{#1 already defined}\fi
        \expandafter\xdef\csname#1\endcsname
        {\the\sectionnumber.\the#2}}%
\def\newsec{
  \global\advance\sectionnumber by 1
  \global\equationnumber=0
  \global\thnumber=0
  \the\sectionnumber .\ }
%
\def\eq#1{\relax
  \global\advance\equationnumber by 1
  \assignnumber{EN#1}\equationnumber
  {\rm \csname EN#1\endcsname}}
\def\eqtag#1{\ifundefined{EN#1}\message{EN#1 undefined}{\sl (#1)}%
  \else\thetag{\csname EN#1\endcsname}\fi}
%
%
\def\thname#1{\relax
  \global\advance\thnumber by 1
  \assignnumber{TH#1}\thnumber
  \csname TH#1\endcsname}
\def\thtag#1{\ifundefined{TH#1}\message{TH#1 undefined}{\sl #1}%
  \else\csname TH#1\endcsname\fi}
%
\comment
\def\eq{}
\def\eqtag#1{(#1)}
\def\thname{}
\def\thtag{}
\endcomment
%
%
\global\countrefno=1

\def\refno#1{\xdef#1{{\the\countrefno}}
\global\advance\countrefno by 1}
\def\R{{\Bbb R}}
\def\N{{\Bbb N}}
\def\C{{\Bbb C}}
\def\CRp{\C\backslash\R_{>0}}

\def\Z{{\Bbb Z}}

\def\Zp{{{\Bbb Z}_{\geq 0}}}

\def\hf{{1\over 2}}

\def\al{\alpha}
\def\be{\beta}

\def\de{\delta}

\def\ep{\epsilon}
\def\si{\sigma}
\def\ta{\tau}
\def\th{\theta}
\def\la{\lambda}
\def\vp{\varphi}

\def\H{\ell^2(\Z)}
\def\Hp{\ell^2(\Zp)}
\def\id{\text{\rm id}}
\refno{\AskeRS}
\refno{\AskeW}
\refno{\GaspCM}
\refno{\GaspR}
\refno{\IsmaS}
\refno{\Kake}
\refno{\KakeMU}
\refno{\KoelSsu}
\refno{\KoelSNATO}
\refno{\KoelVdJSIAM}
\refno{\KoelVdJCA}
\refno{\KoorAM}
\refno{\KoorJF}
\refno{\MasuMNNSU}
\refno{\RoseCM}
\topmatter
\title Transmutation kernels
for the little $q$-Jacobi function transform
\endtitle
\author Erik Koelink and Hjalmar Rosengren
\endauthor
\rightheadtext{transmutation kernels and intertwiners}
\address Technische Universiteit Delft, Faculteit
Informatietechnologie en Systemen, Afd. Toegepaste Wiskundige
Analyse, Postbus 5031, 2600 GA Delft, the Netherlands\endaddress
\email koelink\@dutiaw4.twi.tudelft.nl \endemail
\address Chalmers University of Technology
and G\"oteborg University,
Department of Mathematics, SE-412\, 96 G\"oteborg,
Sweden\endaddress
\email hjalmar\@math.chalmers.se\endemail
\date Version of October 31, 2000
\enddate
\abstract
The little $q$-Jacobi function transform depends on
three parameters. An explicit expression as a sum
of two very-well-poised ${}_8W_7$-series is derived for
the dual transmutation kernel
relating little $q$-Jacobi function transforms
for different parameter sets. A product
formula for the dual transmutation kernel is obtained.
For the inverse transform the transmutation kernel
is given as a ${}_3\vp_2$-series, and a
product formula as a finite sum is derived.
The transmutation kernel gives rise to
intertwining operators for the
second order hypergeometric $q$-difference operator,
which generalise the intertwining operators
arising from a Darboux factorisation.
\endabstract

\endtopmatter

\document

\head \newsec Introduction
\endhead

The Jacobi transform is an integral transform on the
positive half-line with a hypergeometric ${}_2F_1$-series
as its kernel. This transform is a two-parameter
extension of the Fourier-cosine transform and the
Mehler-Fock transform, and also
contains the Hankel transform as a limit case.
The inversion formula for the Jacobi transform can
be found explicitly in several ways,
using asymptotics, spectral analysis, group
theory or intertwining properties.
This transform
has a long history and we refer the reader to the
survey paper \cite{\KoorJF} by Koornwinder.

There are several levels of $q$-analogues of the Jacobi function
and of the corresponding transform pair,
see \cite{\KoelSNATO} for an overview and references.
Here we consider the so-called
little $q$-Jacobi function and the corresponding
transform. The little $q$-Jacobi function transform
has been studied by Kakehi, Masuda and Ueno \cite{\KakeMU},
\cite{\Kake} as the (spherical) Fourier transform on the
quantum $SU(1,1)$ group using the interpretation of
the little $q$-Jacobi functions as matrix elements
of unitary irreducible representations of $U_q(\frak{su}(1,1))$,
see \cite{\MasuMNNSU}. On the other hand the
little $q$-Jacobi function transform occurs when
studying the action of so-called twisted primitive elements
in the principal unitary series representations of
$U_q(\frak{su}(1,1))$, and this has been the motivation
for the study of this paper. However, the paper
is completely analytic in nature, and the quantum group
theoretic interpretation is only discussed briefly in \S 6.

The little $q$-Jacobi function transform can be obtained from
the spectral analysis of the
second order hypergeometric $q$-difference operator
$$
L=L^{(a,b)}= a^2(1+\frac{1}{x})\bigl( T_q-\text{Id}\bigr) +
(1+\frac{aq}{bx})\bigl( T^{-1}_q-\text{Id}\bigl),
\tag\eq{150}
$$
where $T_qf(x)=f(qx)$, on a suitable Hilbert space,
see Kakehi \cite{\Kake} or \cite{\KoelSsu, App.~A}, where
a slightly more general result is given.
So we have eigenfunctions to $L$
in terms of basic hypergeometric series, see
\cite{\GaspR, Ch.~1}.
The little $q$-Jacobi function is defined as
$$
\phi_\la(x;a,b;q) = {}_2\vp_1\left( {{a\si,a/\si}\atop{ab}};
q,-\frac{bx}{a}\right), \quad \la=\hf(\si+\si^{-1}).
\tag\eq{140}
$$
The notation for $q$-hypergeometric series follows
Gasper and Rahman \cite{\GaspR}, and we assume $0<q<1$;
$$
\gather
{}_{k+1}\vp_k\left(
{{a_1,a_2,\ldots,a_{k+1}}\atop{b_1,\ldots,b_k}};
q,z\right) =\sum_{j=0}^\infty
\frac{(a_1,a_2,\ldots,a_{k+1};q)_j}
{(q,b_1,\ldots, b_k;q)_j}z^j,\tag\eq{135}\\
(a_1,\ldots,a_{k+1};q)_j = (a_1;q)_j\ldots (a_{k+1};q)_j,
\ \ (a;q)_j =\prod_{i=0}^{j-1} (1-aq^i),\ j\in\Zp\cup\{\infty\}.
\endgather
$$
The radius of convergence is $1$ for generic parameters, but
there exists a one-valued analytic continuation to
$\C\backslash \R_{\geq 1}$, see \cite{\GaspR, \S 4.5}.

The little $q$-Jacobi function satisfies
$L\phi_\la(\cdot;a,b;q) =
(-1-a^2+2a\la)\phi_\la(\cdot;a,b;q)$.
For later use we note that the little
$q$-Jacobi functions are eigenfunctions for the eigenvalue
$\la$ of
$$
{\Cal L}^{(a,b)} = \frac{1}{2a} L^{(a,b)} + \hf(a+a^{-1})
= \frac{a}{2}(1+\frac{1}{x})T_q
- \bigl( \frac{a}{2x} + \frac{q}{2bx}\bigr)\text{Id}
+ \frac{1}{2a}(1+\frac{aq}{bx})T_q^{-1}.
\tag\eq{155}
$$

For simplicity we assume that $a,b>0$, $ab<1$ and $y>0$,
but the results
hold, mutatis mutandis, for the more general range of the
parameters as discussed in \cite{\KoelSsu, App.~A}. Then
the operator $L$ is an unbounded symmetric operator on the
Hilbert space ${\Cal H}(a,b;y)$ of square integrable
sequences $u=(u_k)_{k\in\Z}$ with respect to the
weights
$$
\sum_{k=-\infty}^\infty |u_k|^2 (ab)^k
\frac{(-byq^k/a;q)_\infty}{(-yq^k;q)_\infty},
\tag\eq{160}
$$
where the operator $L$ is initially defined on the
sequences with finitely many non-zero entries.

Note that \eqtag{160} may be written as a $q$-integral. Indeed,
by associating to $u$
a function $f$ on $yq^\Z$ by $f(yq^k)=u_k$ and using
the notation, see \cite{\GaspR, \S 1.11},
$$
\int_0^{\infty(y)} f(x)\, d_qx = y\sum_{k=-\infty}^\infty
f(yq^k) q^k,
\tag\eq{165}
$$
we see that for
$a=q^{\hf(\al+\be+1)}$ and $b=q^{\hf(\al-\be+1)}$
the sum in \eqtag{160} can be written as
$$
y^\al \int_0^{\infty(y)} f(x) x^\al
\frac{(-xq^{-\be};q)_\infty}
{(-x;q)_\infty}\, d_qx, \qquad \Re \al >-1.
\tag\eq{166}
$$
Using the $q$-binomial theorem, see \cite{\GaspR, \S 1.3},
we see that the quotient of $q$-shifted factorials in
\eqtag{166} tends to $(1+x)^\be$ as $q$ tends to $1$.
In the paper we will use the correspondence
between $u\in{\Cal H}(a,b;y)$
and functions $f$ given by $f(yq^k)=u_k$ repeatedly.

The spectral analysis of $L$, or equivalently ${\Cal L}^{(a,b)}$,
on ${\Cal H}(a,b;y)$
can be completely carried out, and this leads to
corresponding transform
$$
\aligned
\bigl( {\Cal F}_{a,b,y}u\bigr) (\la) &= \sum_{k=-\infty}^\infty
u_k \phi_\la(yq^k;a,b;q) (ab)^k
\frac{(-byq^k/a;q)_\infty}{(-yq^k;q)_\infty}, \\
u_k &= \int_\R  \bigl( {\Cal F}_{a,b,y}u\bigr) (\la)
\phi_\la(yq^k;a,b;q) \, d\nu(\la;a,b;y;q),
\endaligned
\tag\eq{170}
$$
for an explicit measure $d\nu(\cdot;a,b;y;q)$ described
in \thetag{1.10}. 
Here we use the one-valued analytic continuation
of the ${}_2\vp_1$-series.
Then ${\Cal F}={\Cal F}_{a,b,y}$ extends to
an isometric operator from ${\Cal H}(a,b;y)$
onto $L^2\bigl(d\nu(\cdot;a,b;y;q)\bigr)$.

The goal of this paper is to establish a number
of links between two little $q$-Jacobi function transforms
for different parameters $(a,b,y)$.
Our main interest lies with the dual transmutation kernels $P$
satisfying
$$
\bigl({\Cal F}_{c,d,y}[\de_t u]\bigr)(\mu) =
\int_\R  \bigl({\Cal F}_{c,d,y}u\bigr)(\la) \, P_t(\la,\mu)
\, d\nu(\la;a,b;y;q),
\tag\eq{181}
$$
where $(\de_t u)_k=t^ku_k$ for an extra
parameter $t$. Note that
the kernel $P$ is analogous to a non-symmetric Poisson kernel,
and to a Poisson kernel for $(a,b)=(c,d)$, for a family
of orthogonal polynomials. In our main result,
Theorem~2.1, 
we present an explicit expression for the kernel in case
$ad=bc$. The result is much inspired by Mizan Rahman's
summation formulas \cite{\KoelSsu, App.~B}.
We similarly study the transmutation kernel $P$ for the
inverse transform,
$$
\bigl( {\Cal F}^{-1}_{c,d,x}f\bigr)_l =
\sum_{k=-\infty}^\infty \bigl( {\Cal F}^{-1}_{a,b,y}f)_k\,
P_{k,l}\,  (ab)^k
\frac{(-byq^k/a;q)_\infty}{(-yq^k;q)_\infty}.
$$
An explicit expression for the
transmutation kernel is obtained in case $dx=by$
in Theorem~2.2. 
For the transmutation kernels
we show that these kernels do indeed satisfy the
transmutation property, i.e. they intertwine
the second order $q$-difference operator ${\Cal
L}^{(a,b)}$ for different parameters $(a,b)$.
This is closely related to results on $q$-analogues of
Erd\'elyi's fractional integrals recently obtained
by Gasper \cite{\GaspCM}, and we give new proofs of
some of his main results using the little $q$-Jacobi function
transform \eqtag{170}.
The main results are formulated in
\S 2. The proofs of the statements are contained in
\S 3, \S 4 and \S 5, where in \S 3 we prove some summation
formulas, amongst others an extension of
Ramanujan's ${}_1\psi_1$-summation, that are of independent
interest.
In the  final \S 6 we discuss some of the quantum group
theoretic interpretations of these results.

The little $q$-Jacobi functions may also be considered as
$q$-analogues of the Bessel function, see \cite{\KoelSNATO}
for further references. It would be interesting to see if
the (generalised) transmutation kernels can be evaluated also
for other entries in the Askey-Wilson function scheme as
described in \cite{\KoelSNATO}.

We now give the precise form of the
spectral measure $d\nu$ in \eqtag{170}. The measure
can be obtained from the $c$-function expansion for the
little $q$-Jacobi function, see Kakehi \cite{\Kake} or
\cite{\KoelSsu, App.~A}, analogously to the case of the
Jacobi transform. Explicitly,
$$
\aligned
\phi_\la(yq^k;a,b;q) &= c(\si;a,b;q) \Phi_\si(yq^k;a,b;q)
+ c(\si^{-1};a,b;q) \Phi_{\si^{-1}}(yq^k;a,b;q), \\
\Phi_\si(yq^k;a,b;q) &= (a\si)^{-k}
\, {}_2\vp_1\left(
{{a\si,q\si/b}\atop{q\si^2}};q, -\frac{q^{1-k}}{y}\right), \\
c(\si;a,b,y;q) &= \frac{(b/\si,a/\si;q)_\infty}
{(\si^{-2},ab;q)_\infty}\frac{(-by\si, -q/by\si;q)_\infty}
{(-by/a,-qa/by;q)_\infty},
\endgathered
\tag\eq{180}
$$
valid for $\si^2\notin q^\Z$. Then $\Phi_\si$ is the asymptotically
free solution;
$L\Phi_\si(\cdot;a,b;q)=(-1-a^2+2a\la)\Phi_\si(\cdot;a,b;q)$
on $yq^\Z$
with, as before, $\la=\hf(\si+\si^{-1})$. The measure $d\nu$
in \eqtag{170} can be obtained from \eqtag{180}, see \cite{\Kake},
\cite{\KoelSsu, App.~A}. For this we now assume that $a>b$,
which we can do without loss of generality, cf.
Lemma~5.2 
and \thetag{5.8}. 
Explicitly, we have
$$
\multline
\frac{1}{C} \int_\R  f(\la) \, d\nu(\la;a,b;y;q) =
\frac{1}{2\pi}\int_0^\pi f(\cos\th)\, w(e^{i\th})\, d\th + \\
\sum_{k\in\Zp,\, |aq^k|>1}f\Bigl(\hf\bigl(aq^k+a^{-1}q^{-k}\bigr)
\Bigr)\, w_k +
\sum_{k\in\Z,\, |q^{1-k}/by|>1} f\Bigl(-\hf\bigl(\frac{q^{1-k}}{by}+byq^{k-1}\bigr)
\Bigr)\, v_k,
\endmultline
\tag\eq{190}
$$
where
$$
\align
C &=(ab,ab,-by/a,-aq/by,-y,-q/y;q)_\infty, \\
w(z) &= \frac{(z^2,z^{-2};q)_\infty}
{(az,a/z,bz,b/z,-byz, -q/byz,-by/z,-qz/by;q)_\infty}, \\
w_k &= \frac{1-a^2q^{2k}}{1-a^2} \frac{(a^2,ab;q)_k}{(q,aq/b;q)_k}
(ab)^{-k} \frac{(a^{-2};q)_\infty}
{(q,ab,b/a,-aby,-q/aby,-by/a,-aq/by;q)_\infty},\\
v_k &= \frac{\bigl( \frac{q^{2-2k}}{b^2y^2}-1\bigr) q^{-k(k-1)}
(q^2/b^2y^2)^{k-1}}{(q,q,-q^{1-k}/y,-aq^{1-k}/by,
-b^2yq^{k-1},-abyq^{k-1};q)_\infty}.
\endalign
$$
Note that the integral plus the first sum and the sum over $-\Zp$
in the second sum can be written as $dm(\la;a,b,-by,-q/by\mid q)$,
where $dm(\cdot;a,b,c,d|q)$ denotes the standard (non-normalised)
Askey-Wilson measure, see \cite{\AskeW}, \cite{\GaspR, Ch.~6}.
So the measure in \eqtag{190} has an absolutely
continuous part supported on $[-1,1]$, and
for $a<1$ there are no other discrete mass points in $d\nu$
apart from the infinite series tending to $-\infty$,
see \cite{\Kake}, \cite{\KoelSsu, \S 5}, \cite{\KoelSNATO}.
It should be observed that \eqtag{170} can be formally obtained
from the limit transition of the Askey-Wilson polynomials
to the little $q$-Jacobi functions, see
\cite{\KoelSNATO, \S\S 2.3, 4.1, 6.1}.

Later in this paper, especially in \S 3, we
use the notation for very-well-poised series,
cf. \cite{\GaspR, \S 2.1};
$$
\aligned
{}_{r+1}W_r(a_1;a_4,\ldots,a_{r+1};q,z) &=
\sum_{j=0}^\infty  \frac{1-a_1q^{2j}}{1-a_1}
\frac{(a_1,a_4,\ldots,a_{r+1};q)_jz^j }{(q,qa_1/a_4,\ldots,
qa_1/a_{r+1};q)_j}\\ &=
{}_{r+1}\vp_r\left( {{a_1, q\sqrt{a_1},-q\sqrt{a_1},
a_4,\ldots, a_{r+1}}\atop{\sqrt{a_1},-\sqrt{a_1},
qa_1/a_4,\ldots,qa_1/a_{r+1}}};q,z\right).
\endaligned
\tag\eq{137}
$$

\demo{Acknowledgement} The main research for this paper was
done while the second author was employed by Technische
Universiteit Delft.
\enddemo

\head \newsec Statement of main results
\endhead

In this section we describe the main results of the paper.
We start with the dual
transmutation kernel for the little $q$-Jacobi
function, i.e. we want an explicit expression for
$$
\aligned
{\Cal F}_{a,b,y}&\bigl[ k\mapsto
t^k\, \phi_\mu(xq^k;c,d;q)\bigr](\la) \\ &=
\sum_{k=-\infty}^\infty (abt)^k \phi_\mu(xq^k;c,d;q)
\phi_\la(yq^k;a,b;q) \frac{(-byq^k/a;q)_\infty}
{(-yq^k;q)_\infty} \\
&= \sum_{k=-\infty}^\infty (abt)^k
\,{}_2\vp_1\left(
{{c\ta,c/\ta}\atop{cd}};q,-q^k\frac{dx}{c}\right)
\,{}_2\vp_1\left(
{{b\si,b/\si}\atop{ab}};q,-yq^k\right)
\endaligned
\tag\eq{3125}
$$
using \eqtag{170} and Heine's transformation formula
\cite{\GaspR , (1.4.6)}. Here $\la=\hf(\si+\si^{-1})$
and $\mu=\hf(\ta+\ta^{-1})$.
In general we do not have an explicit
expression, but we have the following theorem, which
will be proved in \S 3.

\proclaim{Theorem \thname{210}} Let $a>b>0$, $ab<1$, $y>0$.
Define the dual transmutation kernel
$$
P_t(\la,\mu;q^\al;a,b) ={\Cal F}_{a,b,y}\bigl[ k\mapsto
t^k\, \phi_\mu(yq^k;aq^\al,bq^\al;q)\bigr](\la),
$$
then, using the notation
$\la=\hf(\si+\si^{-1})$ and $\mu=\hf(\ta+\ta^{-1})$ with
$|\si|, |\ta| \geq 1$, the series defining
$P_t$ by \eqtag{170} is absolutely convergent for
$|abq^\al\si\ta|<|abt|<1$.
For $|t|>q/ab$, $t\notin a^2q^{2\al+\Zp}$,
$t\notin b^2q^{2\al+\Zp}$,
$P_t$ can be expressed explicitly as
$$
\align
P_t(\la,\mu;q^\al;a,b) =&
\frac{(b\si^{\pm 1},bq^\al\ta^{\pm 1}, aq^{2\al}\si^{\pm 1}/t,
aq^\al\ta^{\pm 1}/t,q,-yabt,-q/abyt;q)_\infty}
{(ab,abq^{2\al},b/a,q^\al \si^{\pm 1}\ta^{\pm 1}/t, a^2q^{2\al}/t,
abt,-y,-q/y;q)_\infty} \\ &\qquad\qquad\qquad \times
{}_8W_7\left( a^2q^{2\al-1}/t;aq^\al\ta^{\pm 1},a\si^{\pm 1},
abq^{2\al-1}/t;q,q/abt\right) \\
+&
\frac{(a\si^{\pm 1},aq^\al\ta^{\pm 1}, bq^{2\al}\si^{\pm 1}/t,
bq^\al\ta^{\pm 1}/t,q,-b^2yt,-q/b^2yt;q)_\infty}
{(ab,abq^{2\al},a/b,q^\al \si^{\pm 1}\ta^{\pm 1}/t, b^2q^{2\al}/t,
abt,-by/a,-qb/ay;q)_\infty} \\ &\qquad\qquad\qquad \times
{}_8W_7\left( b^2q^{2\al-1}/t;bq^\al\ta^{\pm 1},b\si^{\pm 1},
abq^{2\al-1}/t;q,q/abt\right).
\endalign
$$
The expression for the dual transmutation kernel remains valid for
$\la$ a discrete mass point of the measure $d\nu(\cdot;a,b;y;q)$
as defined in \eqtag{190}.
Moreover, for $q^\al|\ta|<|t|<1/\sqrt{ab}$
and $q^{\be-\al}|t\rho|<|s|< tq^{-2\al}/\sqrt{ab}$,
where $\nu=\hf(\rho+\rho^{-1})$, $|\rho|\geq 1$,
the product formula
$$
P_s(\mu,\nu;q^\be;aq^\al,bq^\al) = \int_\R
P_{q^{2\al}s/t}(\la,\nu;q^{\al+\be};a,b)
P_t(\la,\mu;q^\al;a,b)\, d\nu(\la;a,b;y;q)
$$
is valid.
\endproclaim

As remarked in \S 1
the little $q$-Jacobi functions can be obtained
as a limit of the Askey-Wilson polynomial. In this limit
transition one of the parameters tends to zero and another
one of the parameters
tends to $\infty$ exponentially. This case is not considered
in Askey et al. \cite{\AskeRS, \S 3}, where the non-symmetric
Poisson kernel for the Askey-Wilson polynomials is derived.
Note that the expression \cite{\AskeRS, (3.9)-(3.11)} is
much more complicated than the expression in Theorem \thtag{210}.
Motivated by the quantum group theoretic interpretation,
see \S 6,
Theorem \thtag{210} should be compared to the non-symmetric
Poisson kernel for Al-Salam and Chihara polynomials, which
consists of one very-well-poised ${}_8W_7$-series,
see Askey, Rahman and Suslov \cite{\AskeRS, (14.8)}
and Ismail and Stanton \cite{\IsmaS, Thm.~4.2}.

For the transmutation kernel we have the following result.

\proclaim{Theorem \thname{230}} Define the transmutation
kernel
$$
\multline
P_{k,l}(a,b,y;r,s) =
{\Cal F}^{-1}_{a,b,y}[ \la\mapsto \phi_\la(yq^l/s;ar,bs;q)]_k
\\ = \int_\R
\phi_\la(yq^l/s;ar,bs;q)\phi_\la(yq^k;a,b;q)\,
d\nu(\la;a,b;y;q).
\endmultline
$$
For $r,s>0$, $rs<1$,
the transmutation kernel is given by
$$
P_{k,l}(a,b,y;r,s)
= (ab)^{-l}
\frac{(ab,rs,q^{k-l+1},-yq^k;q)_\infty}
{(q,abrs,rsq^{k-l},-byq^k/ar;q)_\infty}
\, {}_3\vp_2\left( {{q^{l-k},\, r,\, ar/b}
\atop{rs, -arq^{1-k}/by}};q,q\right)
$$
with the convention $P_{k,l}(a,b,y;r,s)=0$
for $k<l$. Moreover, the transmutation kernel
satisfies the product formula for $r,s,t,u>0$, $rs<1$,
$tu<1$,
$$
\sum_{l=k}^p P_{k,l}(a,b,y;r,s)
P_{l,p}(ar,bs,\frac{y}{s};rt,su)\, (abrs)^l
\frac{(-ybq^l/ar;q)_\infty}
{(-yq^l/s;q)_\infty} = P_{k,p}(a,b,y;rt,su).
$$
\endproclaim

See \thetag{4.7} 
for the explicit expression for the
product formula in terms of the ${}_3\vp_2$-series.
The resulting product formula is already contained
in Gasper \cite{\GaspCM, (1.7)}, and it can also be
obtained from a general expansion formula
\cite{\GaspR, (3.7.9) with $k=r=t=u=2$, $s=1$} due to Verma
using transformation formulas for ${}_3\vp_2$-series,
see also \cite{\GaspCM, \S 4}.

Note that for $r=1$ (or $s=1$) the
transmutation kernel in Theorem \thtag{230} simplifies;
for $k\geq l$ and $|s|<1$
$$
P_{k,l}(a,b,y;1,s)
= (ab)^{-l} \frac{(ab,-yq^k;q)_\infty}
{(abs,-byq^k/a;q)_\infty} \frac{(s;q)_{k-l}}{(q;q)_{k-l}}
$$
and letting $s\uparrow 1$ gives
$$
\lim_{s\uparrow 1} P_{k,l}(a,b,y;1,s) =
\de_{k,l} (ab)^{-k} \frac{(-yq^k;q)_\infty}
{(-byq^k/a;q)_\infty}
$$
in accordance with \eqtag{160} and \eqtag{170}.

The last main result deals with intertwining operators
for the second order $q$-difference operator ${\Cal L}^{(a,b)}$
as in \eqtag{155}. The intertwining properties are also
known as transmutation properties and for the Jacobi function
transform the intertwining operators are known as the
Abel transform, see \cite{\KoorJF}.

\proclaim{Theorem \thname{240}} {\rm (i)} Let
$a,b\in\C\backslash\{0\}$, $\nu,\mu\in\C$ with
$|q^{\nu-\mu}b/a|<1$. Define the operator
$$
\multline
\bigl(W_{\nu,\mu}(a,b)f\bigr)(x) =
\frac{(-x;q)_\infty}{(-xq^{-\mu};q)_\infty}
q^{-\mu^2}\bigl( \frac{b}{a}\bigr)^\mu x^{\mu+\nu}
\\ \times \sum_{p=0}^\infty
f(xq^{-\mu-p})\, q^{-p\nu}\frac{(q^\nu;q)_p}{(q;q)_p}
\, {}_3\vp_2\left( {{q^{-p}, q^{-\mu},-q^{1+\mu-\nu}a/bx}\atop
{q^{1-p-\nu}.-q^{\mu+1}/x}};q,q^{1-\mu}\frac{b}{a}\right)
\endmultline
$$
for any function $f$ with $|f(xq^{-p})|={\Cal O}(q^{p(\ep+\nu)})$
for some $\ep>0$.
Then $W_{\nu,\mu}(a,b)\circ
{\Cal L}^{(a,b)}= {\Cal L}^{(aq^{-\nu},bq^{-\mu})}
\circ W_{\nu,\mu}(a,b)$ on the space of compactly supported
functions and for $|a\si|<q^\nu$
$$
\bigl(W_{\nu,\mu}(a,b)\Phi_\si(\cdot;a,b;q)\bigr)(yq^k) =
y^{\mu+\nu}
\frac{(a\si,b\si;q)_\infty}
{(aq^{-\nu}\si,bq^{-\mu}\si;q)_\infty}
\Phi_\si(yq^k;aq^{-\nu},bq^{-\mu};q).
$$
{\rm (ii)}
Let $a,b>0$, $ab<1$, $\nu>0$ and
$\mu\in\C\backslash\Z_{\leq 0}$. Define the operator
$$
\multline
\bigl( A_{\nu,\mu}(a,b)f\bigr)(x) =
\frac{(-bxq^\mu/a;q)_\infty}{(-bxq^{\mu-\nu}/a;q)_\infty}
\\ \times \sum_{k=0}^\infty f(xq^{\mu+k})\,
(ab)^k \frac{(q^\nu,-xq^\mu;q)_k}{(q,-bxq^\mu/;q)_k}
\, {}_3\vp_2\left( {{q^{-k}, q^\mu,-bxq^{\mu-\nu}/a}
\atop{q^{1-\nu-k},-xq^\mu}};q,q\right)
\endmultline
$$
for any bounded function. Then
${\Cal L}^{(aq^\nu,bq^\mu)}\circ A_{\nu,\mu}(a,b) =
A_{\nu,\mu}(a,b)\circ {\Cal L}^{(a,b)}$
on the space of functions compactly supported in
$(0,\infty)$. Moreover,
$$
\bigl(A_{\nu,\mu}(a,b) \phi_\la(\cdot;a,b;q)\bigr)(x) =
\frac{(abq^{\nu+\mu};q)_\infty}{(ab;q)_\infty}
\, \phi_\la(x;aq^\nu,bq^\mu;q).
\tag\eq{250}
$$
\endproclaim

The operators $W_{\nu,\mu}(a,b)$ and $A_{\nu,\mu}(a,b)$
are $q$-analogues of the (generalised) Abel transform,
see \cite{\KoorJF, \S 5}. Note that Theorem \thtag{240}
gives $q$-integral representations for the little $q$-Jacobi
function and the asymptotically free solution $\Phi_\si$
of \eqtag{180}.

The ${}_3\vp_2$-kernel of $A_{\nu,\mu}(a,b)$
is the same as the transmutation kernel.
In order to see this we
first invert the summation for the ${}_3\vp_2$-series in
\eqtag{250} using \cite{\GaspR, exerc. 1.4(ii)} and next transform
it using \cite{\GaspR, (III.13)}. Hence
\eqtag{250} is equivalent to
\thetag{4.1}. 
This shows that the transform with the transmutation kernel
does indeed satisfy the transmutation property.

\head \newsec The dual transmutation kernel
\endhead

The goal of this section is to prove Theorem \thtag{210}.
We start with proving some general results, which are
of independent interest, and come back to the proof of
Theorem \thtag{210} later.

We first formulate a general proposition generalising
Rahman's summation formulas in \cite{\KoelSsu, App.~B}.
Note that in case the argument of the basic hypergeometric
series has absolute value bigger than $1$, we implicitly
use the one-valued analytic continuation to
$\C\backslash\R_{\geq1}$.

\proclaim{Proposition \thname{310}}  Let $x,y\in \CRp$.
Consider the sum
$$
S=\sum_{n=-\infty}^\infty z^n {}_2\vp_1\left( {{a,b}\atop{c}};q,
xq^n\right) {}_2\vp_1\left( {{d,e}\atop{f}};q, yq^n\right),
$$
which is absolutely convergent for
$\max(|ad|,|ae|,|bd|,|be|)<|z|<1$. If furthermore
$abde=cf$, $fx=dey$, $q<|z|$, and $z/abf, z/cde\not\in q^\Zp$
then $S$ equals
$$
\multline
\frac{(e,d,c/a,c/b,abd/z,abe/z,af/z,bf/z,q,yz,q/yz;q)_\infty}
{(c,f,c/ab,ae/z,be/z,ad/z,bd/z,abf/z,z,y,q/y;q)_\infty}
\, {}_8W_7\bigl( \frac{abf}{qz};a,b,\frac{f}{e},\frac{f}{d},
\frac{cf}{qz};q, \frac{q}{z}\bigr) \\
+
\frac{(a,b,f/d,f/e,ade/z,bde/z,cd/z,ce/z,q,xz,q/xz;q)_\infty}
{(c,f,f/de, ae/z,be/z,ad/z,bd/z,
cde/z,z,x,q/x;q)_\infty}
\, {}_8W_7\bigl( \frac{cde}{qz};d,e,\frac{c}{a},\frac{c}{b},
\frac{cf}{qz};q, \frac{q}{z}\bigr) .
\endmultline
$$
\endproclaim

\demo{Remark \thname{311}}
Using the transformation formulas for very-well-poised
series we can find more 
expressions for the sum $S$, some of
which are valid also if one of the conditions
$q<|z|$, $z/abf, z/cde\notin q^\Zp$ is violated,
see \thetag{3.11}, 
\thetag{3.12},  
\thetag{3.13} below. 
Here we have chosen
the expression as a sum of two very-well-poised series
that shows the symmetries $(a,b,c,x)\leftrightarrow
(d,e,f,y)$, $a\leftrightarrow b$,
and $d\leftrightarrow e$, which are
obvious in the sum. (Note that the conditions $abde=cf$, $fx=dey$
also display this symmetry.) The expression also
displays the symmetry $x\leftrightarrow y$,
$(a,b,d,e) \leftrightarrow (c/a,c/b,f/e,f/d)$.
For the right hand sides this follows from $abde=cf$ and $fx=dey$.
For the sum $S$ this follows
from a double application of
Heine's transformation \cite{\GaspR, (III.3)}, since
$fx=dey$ and $abx=cy$.
\enddemo

The following lemma is of use in the proof of Proposition
\thtag{310}, and is of independent interest. Note that in
case $k=0$ the series in the summand can be
summed by the $q$-binomial formula, and we obtain
Ramanujan's ${}_1\psi_1$-summation formula, see \cite{\GaspR,
(5.2.1)}.

\proclaim{Lemma \thname{320}}
For $\max(|a_1|,\ldots,|a_{k+1}|)<|z|<1$ and for
$x\in\CRp$ we have
$$
\sum_{n=-\infty}^\infty
 {}_{k+1}\vp_k\left(
{{a_1,\ldots,a_{k+1}}\atop{b_1,\ldots,b_k}};q,xq^n\right) z^n=
\frac{(a_1,\ldots,a_{k+1},b_1/z,\ldots,b_k/z,q,xz,q/xz;q)_\infty}
{(b_1,\ldots,b_k,a_1/z,\ldots,a_{k+1}/z,z,x,q/x;q)_\infty}.
$$
\endproclaim

\demo{Proof} By shifting the summation parameter
we can assume, without loss of generality, that $1\leq |x|<q^{-1}$.
The series $\sum_{n=1}^\infty$ is absolutely convergent
for $|z|<1$. Let us assume for the moment that
$|qb_1\ldots b_k|<|xq^na_1\ldots a_{k+1}|$ for $n\leq0$,
then the analytic
continuation of the ${}_{r+1}\vp_r$-series in the summand
for $|xq^n|\geq 1$ is given by
$$
\aligned
{}_{k+1}\vp_k&\left(
{{a_1,\ldots,a_{k+1}}\atop{b_1,\ldots,b_k}};q,xq^n\right) = \\
&\frac{(a_2,\ldots,a_{k+1},b_1/a_1,\ldots,b_k/a_1,
a_1xq^n,q^{1-n}/a_1x;q)_\infty}{(b_1,\ldots,b_k,a_2/a_1,\ldots,
a_{k+1}/a_1, xq^n,q^{1-n}/x;q)_\infty} \\
&\times {}_{k+1}\vp_k \left( {{a_1,qa_1/b_1,\ldots,qa_1/b_k}
\atop{qa_1/a_2,\ldots,qa_1/a_{k+1}}};q,
\frac{q^{1-n}b_1\ldots b_k}{xa_1\ldots a_{k+1}}\right) \\
&+  \text{\rm idem}(a_1;a_2,\ldots,a_{k+1}),
\endaligned
\tag\eq{330}
$$
where $\text{\rm idem}(a_1;a_2,\ldots,a_{k+1})$ after an expression
stands for the sum of $k$ terms obtained from the previous
expression by interchanging $a_1$ with each $a_i$, $i=2,3,\ldots,
k+1$, see \cite{\GaspR, (4.5.2)}. Note that \eqtag{180} is the
case $k=1$ of \eqtag{330}.

Using the theta product identity
$$
(aq^n,q^{1-n}/a;q)_\infty = (-a)^{-n} q^{-\hf n(n-1)}
(a,q/a;q)_\infty
\tag\eq{340}
$$
we see that the $n$-dependence in \eqtag{330} simplifies.
Indeed, since
$$
\frac{(a_1xq^n,q^{1-n}/a_1x;q)_\infty}
{(xq^n,q^{1-n}/x;q)_\infty} = a_1^{-n}
\frac{(a_1x,q/a_1x;q)_\infty}
{(x,q/x;q)_\infty}
\tag\eq{350}
$$
the sum $\sum_{n=-\infty}^0$ of the first term in
the right hand side of \eqtag{330}
times $z^n$ is absolutely convergent for $|z/a_1|>1$. Hence,
the sum is absolutely convergent for
$\max(|a_1|,\ldots,|a_{k+1}|)<|z|<1$.

Next we split the sum $\sum_{n=1}^\infty +\sum_{n=-\infty}^0$,
using the series for the ${}_{k+1}\vp_k$-series in the first
sum and \eqtag{330} in the second sum. Interchanging summations
we see from \eqtag{350} that the sums over $n$ are all geometric.
So we see that the left hand side of the lemma equals
$$
\aligned
&z\sum_{j=0}^\infty \frac{(a_1,\ldots,a_{k+1};q)_j}
{(q,b_1,\ldots,b_k;q)_j} \frac{(qx)^j}{1-q^jz}\\
+&\frac{(a_2,\ldots,a_{k+1},b_1/a_1,\ldots,b_k/a_1,
a_1x,q/a_1x;q)_\infty}{(b_1,\ldots,b_k,a_2/a_1,\ldots,
a_{k+1}/a_1, x,q/x;q)_\infty}\\
&\times \sum_{j=0}^\infty
\frac{(a_1,qa_1/b_1,\ldots,qa_1/b_k;q)_j}{(q,qa_1/a_2,\ldots,
qa_1/a_{k+1};q)_j} \left( \frac{qb_1\ldots b_k}{xa_1\ldots a_{k+1}}
\right)^j \frac{1}{1-a_1q^j/z} \\
+&\text{\rm idem}(a_1;a_2,\ldots,a_{k+1}).
\endaligned
\tag\eq{360}
$$
The sums in \eqtag{360}
can be written as ${}_{k+2}\vp_{k+1}$-series.
The first sum in \eqtag{360} equals
$$
\aligned
\frac{z}{z-1}& \, {}_{k+2}\vp_{k+1}\left( {{a_1,\ldots, a_{k+1},z}
\atop{b_1,\ldots, b_k,qz}};q,qx\right) = \\
&\frac{z}{z-1} \frac{(a_1,\ldots,a_{k+1},b_1/z,\ldots,b_k/z,
q, xqz,1/xz;q)_\infty}{(b_1,\ldots,b_k,a_1/z,\ldots,
a_{k+1}/z,qz, xq,1/x;q)_\infty} \\
+ & \frac{z}{z-1}
\frac{(a_2,\ldots,a_{k+1},b_1/a_1,\ldots,b_k/a_1,
a_1xq,1/a_1x,z,qz/a_1;q)_\infty}{(b_1,\ldots,b_k,a_2/a_1,\ldots,
a_{k+1}/a_1, xq,1/x, qz,z/a_1;q)_\infty} \\
&\times {}_{k+2}\vp_{k+1} \left(
{{a_1,qa_1/b_1,\ldots,qa_1/b_k,a_1/z}
\atop{qa_1/a_2,\ldots,qa_1/a_{k+1},qa_1/z}};q,
\frac{qb_1\ldots b_k}{xa_1\ldots a_{k+1}}\right) \\
+ & \text{\rm idem}(a_1;a_2,\ldots,a_{k+1}),
\endaligned
\tag\eq{370}
$$
where we have used \cite{\GaspR, (4.5.2)}, cf. \eqtag{330},
once again. (The first ${}_{k+2}\vp_{k+1}$-series reduces to $1$,
since an upper parameter is equal to $1$.) The $k+1$
${}_{k+2}\vp_{k+1}$-series in \eqtag{370}
are the same as in \eqtag{360}, and a simple calculation
reveals that they occur with
opposite coefficients. Hence, using \eqtag{370} in \eqtag{360}
leaves only the first term on the right hand side
of \eqtag{370}. This proves the result for the condition
$|qb_1\ldots b_k|<|xq^na_1\ldots a_{k+1}|$ for $n\leq0$,
and the general case follows by analytic continuation in the
parameters of the ${}_{k+1}\vp_k$-series.
\qed\enddemo

\demo{Proof of Proposition \thtag{310}}
The sum $\sum_{n=0}^\infty$ in $S$ is absolutely convergent for
$|z|<1$. As in the proof of Lemma \thtag{320} we can use \eqtag{330}
for $k=1$ twice to see that the sum $\sum_{n=-\infty}^{-1}$ in
$S$ is absolutely convergent for
$|z|>\max (|ad|, |ae|, |bd|, |be|)$.
For $n$ large enough we have $|xq^n|, |yq^n|<1$, so let us
assume first $|x|<1$ and $|y|<1$. We can use
the series representation \eqtag{135} twice to write
$$
{}_2\vp_1\left( {{a,b}\atop{c}};q,
x\right) {}_2\vp_1\left( {{d,e}\atop{f}};q, y\right) = 
\sum_{k=0}^\infty  \frac{(d,e;q)_k}{(q,f;q)_k} y^k
\, {}_4\vp_3\left( {{a,b,q^{-k},q^{1-k}/f}\atop{c,
q^{1-k}/d,q^{1-k}/e}};q, \frac{qfx}{dey}\right),
\tag\eq{380}
$$
cf. Rahman's proof in \cite{\KoelSsu, App.~B}. The terminating
${}_4\vp_3$-series in the summand in \eqtag{380} is balanced
for $abde=cf$ and $dey=fx$, the assumptions in the proposition.
Hence it can be transformed into a terminating very-well-poised
${}_8W_7$-series, see \eqtag{137} for the notation, by
\cite{\GaspR, (III.19)}. It follows that \eqtag{380}
equals
$$
\sum_{k=0}^\infty  \frac{(e,ab,bd;q)_k}{(q,f,abd;q)_k}
y^k \, {}_8W_7\bigl( \frac{abd}{q};a,b,\frac{f}{e},
q^{-k},q^{k-1}cf;q, \frac{qd}{f}\bigr).
\tag\eq{390}
$$

The ${}_8W_7$-series in \eqtag{390} can be rewritten as a sum
of two non-terminating ${}_8W_7$-series using Bailey's
three-term transformation \cite{\GaspR, (III.37)
with $(a,b,c,d,e,f)$ replaced by $(afq^k/e,adq^k,aq/c,q/e,a,f/e)$}.
Recalling that $abde=cf$ we find
$$
\multline
\frac{(afq^k,q^{k+1}f/e,aq^{k+1},eq^k,c/a,c/b,d;q)_\infty}
{(afq^{k+1}/e;q)_\infty}
\, {}_8W_7\bigl( \frac{afq^k}{e};adq^k,\frac{aq}{c},\frac{q}{e},
a,\frac{f}{e};q,beq^k\bigr) = \\
\frac{(fq^k,aeq^k,q^{k+1},c,bd,c/ab,ad;q)_\infty}
{(abd;q)_\infty}\, {}_8W_7\bigl( \frac{abd}{q};a,b,\frac{f}{e},
q^{-k},q^{k-1}cf;q, \frac{qd}{f}\bigr) \\
+ \frac{ed}{f} \frac{(edq/f,cdq^k,adeq^k,dq^{k+1},
cq^{k+1}/b,q/e,a,f/e,b,dq^{1-k}/f,q^kf/d;q)_\infty}
{(qf/de,fcq^k/e,qd/f,q^{1-k}/e,cdq^{k+1}/b;q)_\infty} \\ \times
\, {}_8W_7\bigl( \frac{dcq^k}{b};adq^k,\frac{dq}{f},\frac{q}{b},
d,\frac{c}{b};q,beq^k\bigr).
\endmultline
\tag\eq{3100}
$$
Note that the other two ${}_8W_7$-series in \eqtag{3100}
are obtained from each other by interchanging
$(a,b,c,x)$ with $(d,e,f,y)$. Using \eqtag{3100} in \eqtag{390}
and recalling that $abde=cf$ and $dey=fx$
leads to the following expression for \eqtag{380};
$$
\multline
\sum_{k=0}^\infty y^k
\frac{(abdq^k,afq^k,q^{k+1}f/e,aq^{k+1},
e,c/a,c/b,d;q)_\infty}
{(f,aeq^k,q,c,bdq^k,c/ab,adq^k,afq^{k+1}/e;q)_\infty}
\\ \times
\, {}_8W_7\bigl( \frac{afq^k}{e};adq^k,\frac{aq}{c},\frac{q}{e},
a,\frac{f}{e};q,beq^k\bigr)
+  {\text{\rm Idem}}\bigl( (a,b,c,x);(d,e,f,y)\bigr)
\endmultline
\tag\eq{3110}
$$
where ${\text{\rm Idem}}\bigl( (a,b,c,x);(d,e,f,y)\bigr)$ 
means
that we have the same sum with the parameter sets 
$(a,b,c,x)$ and
$(d,e,f,y)$ interchanged.

Since we assume $|y|<1$ and
$|be|<1$ the double series in \eqtag{3110} is absolutely
convergent. Interchanging summations and
recalling $abde=cf$ shows that
$$
\multline
{}_2\vp_1\left( {{a,b}\atop{c}};q,
x\right) {}_2\vp_1\left( {{d,e}\atop{f}};q, y\right) = \\
\sum_{j=0}^\infty (be)^j 
\frac{(e,c/a,c/b,d,afq^j, q^{j+1}f/e,
aq^{j+1},abdq^j;q)_\infty}
{(f,q,c,c/ab,ae,bd,afq^j/e,adq^j;q)_\infty}
\frac{(aq/c,q/e,a,f/e;q)_j}{(q,fq/de;q)_\infty} \,
\\ \times(1-\frac{afq^{2j}}{e})\
{}_6\vp_5\left( {{q,ae,bd,afq^j/e,adq^j,afq^{1+2j}/e}\atop
{abdq^j,afq^j,q^{j+1}f/e,aq^{j+1},afq^{2j}/e}};q,yq^j\right)\\
+  {\text{\rm Idem}}\bigl( (a,b,c,x);(d,e,f,y)\bigr)
\endmultline
\tag\eq{3120}
$$
assuming that $abde=cf$, $dey=fx$, $|x|<1$, $|y|<1$,
and $|be|<1$. As a function of $y$ the left hand side has a
unique analytic continuation to $\C\backslash(\R_{\geq 1}
\cup \frac{de}{f}\R_{\geq 1})$. The
$y$-dependence in the first sum on the
right hand side is only at the argument spot
of the ${}_6\vp_5$-series, 
which has a unique analytic continuation
to $\C\backslash\R_{\geq 1}$. As $j\to\infty$ the
${}_6\vp_5$-series tends to $1$, so the 
convergence with respect to
$y$ in the right hand side of \eqtag{3120} 
is uniform on compact
sets. Similarly, the ${}_6\vp_5$-series in the other sum
has a unique analytic extension to $\C\backslash
(\frac{de}{f}\R_{\geq 1})$.
Hence \eqtag{3120} remains valid for 
$y\in \C\backslash(\R_{\geq 1}
\cup \frac{de}{f}\R_{\geq 1})$
after using the analytic continuation
of the ${}_2\vp_1$- and ${}_6\vp_5$-series. 
In particular, this
means that \eqtag{3120} is valid for $abde=cf$, $dey=fx$,
$|be|<1$ and  $y\in \C\backslash(\R_{\geq 1}
\cup \frac{de}{f}\R_{\geq 1})$.

To prove the result we replace $x$ and $y$ by $xq^n$ and
$yq^n$ in \eqtag{3120}, multiply by $z^n$ and
sum over $n\in\Z$. If we assume for the moment that
$\max( q, |ae|, |bd|, |ad|, |af/e|, |cd/b|)<|z|<1$, then we
can interchange summations and use Lemma \thtag{320} twice
to sum the inner sum. Some cancellation occurs and after using
the theta product identity \eqtag{340} twice we see that $S$ equals
$$
\aligned
&\frac{(e,c/a,c/b,d,yz,q/yz,q;q)_\infty}
{(f,c,c/ab,z,q/z,ae/z,bd/z,y,q/y;q)_\infty} \\ \times&
\sum_{j=0}^\infty
\Bigl( \frac{be}{z}\Bigr)^j (1-\frac{afq^{2j}}{ez})
\frac{(aq/c,q/e,a,f/e;q)_j}{(q,fq/de;q)_j}
\frac{(abdq^j/z,afq^j/z,fq^{j+1}/ez,aq^{j+1}/z;q)_\infty}
{(afq^j/ez,adq^j/z;q)_\infty}\\
+  &{\text{\rm Idem}}\bigl( (a,b,c,x);(d,e,f,y)\bigr)\\
=&\frac{(e,c/a,c/b,d,q,yz,q/yz,abd/z,af/z,qf/ez,aq/z;q)_\infty}
{(f,c,c/ab,y,q/y,z,q/z,ae/z,bd/z,qaf/ze,ad/z;q)_\infty}
{}_8W_7\bigl( \frac{af}{ze}; \frac{aq}{c},\frac{q}{e},a,\frac{f}{e},
\frac{ad}{z};q, \frac{be}{z}\bigr) \\
+  &{\text{\rm Idem}}\bigl( (a,b,c,x);(d,e,f,y)\bigr).
\endaligned
\tag\eq{3121}
$$
For $q<|z|$ we can use \cite{\GaspR, (III.23)}
using $abde=cf$
to transform the ${}_8W_7$-series in the required form.
Finally, use analytic continuation in $z$ to find the result.
\qed\enddemo

\demo{Remark \thname{3122}}
{\rm (i)} A proposition of this type has been proved first
by Mizan Rahman for two special cases, see \cite{\KoelSsu, App.~B}.
To see how the two special cases are contained in the
result we consider it in the form \eqtag{3121}.
To the first ${}_8W_7$-series we apply \cite{\GaspR, (III.37) with
$(a,b,c,d,e,f)$ replaced by
$(af/ze,ad/z,aq/c,f/e,a,q/e)$} to write it as a sum of two
very-well-poised ${}_8W_7$-series, of which the second is
the same ${}_8W_7$-series as the second in \eqtag{3121}. It turns
out that we can add the coefficients using the theta product
identity as in \cite{\GaspR, Exerc. 2.16
with $(x,\la,\mu,\nu)$ replaced by $(\sqrt{xyz},
\sqrt{xz/y},\sqrt{ef/dz},\sqrt{x/yz})$}.
This yields
$$
\multline
S = \frac{(q,d,e,abd/z,f/z,adq/f,bdq/f,yz,q/yz;q)_\infty}
{(f,z,ad/z,bd/z,e/z,qd/f,qc/e,y,q/y;q)_\infty}
\,{}_8W_7(\frac{abd}{f};
a,b,\frac{q}{e},\frac{c}{z},\frac{qz}{f};q,d) \\
+ \frac{(q,a,b,f/e,qde/f,qd/z,cd/z,ade/z,cq/bz,ey,q/ey,
xz/e,qe/xz;q)_\infty}{(c,f,ad/z,bd/z,ae/z,qd/f,cdq/bz,
y,q/y,x,q/x,z/e,qe/z;q)_\infty}  \\ \times
\, {}_8W_7(\frac{cd}{bz};d,\frac{dq}{f},\frac{q}{b},\frac{c}{b},
\frac{ad}{z};q,\frac{be}{z})
\endmultline
\tag\eq{3123}
$$
valid for $|d|<1$ and on any subregion of
$\max(|ad|,|ae|,|bd|,|be|)<|z|<1$
in the complex $z$-plane as long as the right hand side
is analytic in this subregion.
The first case proved by Mizan Rahman \cite{\KoelSsu, App.~B.1}
corresponds to $f=qz$ in \eqtag{3123},
so that the first ${}_8W_7$-series
reduces to $1$, and the second case proved by Mizan Rahman
\cite{\KoelSsu, App.~B.3} corresponds to $ey\in q^\Z$
in \eqtag{3123}, so that
the second term vanishes.
\par
{\rm (ii)}
The proof of Proposition \thtag{310} is inspired by
Rahman's method as presented in \cite{\KoelSsu, \S B.3},
but is of a different nature.
In Rahman's case the balanced ${}_4\vp_3$-series is written
as a $q$-integral using \cite{\GaspR, (2.10.19)}, which is
a three-term transformation for balanced ${}_4\vp_3$-series.
In this paper we use Bailey's three-term transformation
\cite{\GaspR, (2.11.1)} for
very-well-poised ${}_8\vp_7$-series, which can be deduced from
\cite{\GaspR, (2.10.19)}. Note that the result
as a sum of two very-well-poised series cannot be simplified.
\par
{\rm (iii)} Note that \eqtag{3121} and \eqtag{3123} give
alternative expressions for the sum $S$ of Proposition
\thtag{310}, and together with the obvious symmetries in
$S$ we find more expressions for $S$ in terms of a sum of two
very-well-poised ${}_8W_7$-series. We can also rewrite the
result as a sum of three balanced ${}_4\vp_3$-series as follows.
Start with \eqtag{3123} with $(a,b,c,x)\leftrightarrow (d,e,f,y)$,
apply \cite{\GaspR, (III.36) with $(a,b,c,d,e,f)$
replaced by $(ade/c,q/b,f/z,d,e,qz/c)$} to the first
${}_8W_7$-series to write it as a sum of two balanced
${}_4\vp_3$-series
and apply \cite{\GaspR, (III.36) with $(a,b,c,d,e,f)$
replaced by $(af/ez,q/e,f/e,aq/c,a,ad/z)$} to the second
${}_8W_7$-series to write it as a sum of two balanced
${}_4\vp_3$-series of which one also occurs in the previous
transformation. The balanced ${}_4\vp_3$-series can be taken
together using the theta product identity \cite{\GaspR, Exerc. 2.16
with $(x,\la,\mu,\nu)$ replaced by $(\sqrt{bxz/a},\sqrt{axz/b},
\sqrt{abx/z},\sqrt{yz/c})$} resulting in the following
expression:
$$
\aligned
S =&  \frac{(q,a,b,c/z,xz,q/xz;q)_\infty}
{(c,z,a/z,b/z,x,q/x;q)_\infty} \, {}_4\vp_3\left(
{{d,e,z,qz/c}\atop{f,qz/b,qz/a}};q,q\right) \\
&+ \frac{(q,d,b,c/a,af/z,e,ax,q/ax,yz/a,qa/yz;q)_\infty}
{(c,f,b/a,z/a,ad/z,ae/z,x,q/x,y,q/y;q)_\infty}
\, {}_4\vp_3\left(
{{a,qa/c,ad/z,ae/z}\atop{qa/b,qa/z,af/z}};q,q\right) \\
 &+ \frac{(q,d,a,c/b,bf/z,e,bx,q/bx,yz/b,qb/yz;q)_\infty}
{(c,f,a/b,z/b,bd/z,be/z,x,q/x,y,q/y;q)_\infty}
\, {}_4\vp_3\left(
{{b,qb/c,bd/z,be/z}\atop{qb/a,qb/z,bf/z}};q,q\right),
\endaligned
\tag\eq{3124}
$$
valid on any subregion of $\max(|ad|,|ae|,|bd|,|be|)<|z|<1$
in the complex $z$-plane as long as the right hand side
is analytic in this subregion.
\enddemo

As observed the case $ey\in q^\Z$ is special, and we see
later that this case corresponds to the infinite set of discrete
mass points in the spectral measure $d\nu$ of the little
$q$-Jacobi function transform. However, $ey=q^{1-l}$ with
$l\to\infty$ violates the conditions for absolute convergence
of $S$ as given in Proposition \thtag{310}, so we have to deal
with this case separately. Using Heine's transformation
\cite{\GaspR, (III.1)} we see in this case
$$
\aligned
{}_2\vp_1\left( {{d,e}\atop{f}};q, yq^n\right) &=
\frac{(d,eyq^n;q)_\infty}{(f,yq^n;q)_\infty}
{}_2\vp_1\left( {{f/d,yq^n}\atop{eyq^n}};q, d\right) \\ & =
\frac{(d,q^{1-l+n};q)_\infty}{(f,yq^n;q)_\infty}
{}_2\vp_1\left( {{f/d,yq^n}\atop{q^{1-l+n}}};q, d\right),
\endaligned
\tag\eq{3126}
$$
initially for $|yq^n|<1$ and by analytic continuation in $y$
to the general case for $n\in\Z$ and $y\in\CRp$.
Note that the right hand side of \eqtag{3126} displays
$q$-Bessel coefficient behaviour, which means that for $l>n$
the series on the right hand side of \eqtag{3126} starts
at $l-n$, see the proof of
Proposition~3.5. 

We now prove the necessary
result in some greater generality.

\proclaim{Proposition \thname{331}} Let
$\max (|ya|, |yb|)<|z|<1$, $|y|<1$, $x\in\CRp$
and $|qcd|<|xab|$, then
$$
\align
\sum_{k=-\infty}^\infty &z^k \, {}_2\vp_1\left(
{{a,b}\atop{c}};q,xq^k\right)
\frac{(q^{k+1};q)_\infty}{(dq^k;q)_\infty}
\, {}_2\vp_1\left( {{dq^k,e}\atop{q^{k+1}}};q,y\right) = \\
&\frac{(q,b,c/a,dz/a,ax,q/ax,eya/z;q)_\infty}
{(d,c,b/a,z/a,x,q/x,ya/z;q)_\infty}
\, {}_4\vp_3\left( {{a,aq/c,aq/dz,ya/z}\atop{aq/b,aq/z,eya/z}};
q, \frac{qcd}{xab}\right) \\
+&\frac{(q,a,c/b,dz/b,bx,q/bx,eyb/z;q)_\infty}
{(d,c,a/b,z/b,x,q/x,yb/z;q)_\infty}
\, {}_4\vp_3\left( {{b,bq/c,bq/dz,yb/z}\atop{bq/a,bq/z,eyb/z}};
q, \frac{qcd}{xab}\right) \\
+& \frac{(q,a,b,c/z,zx,q/zx,ey;q)_\infty}
{(z,c,a/z,b/z,x,q/x,y;q)_\infty}
\, {}_4\vp_3\left( {{z,zq/c,q/d,y}\atop{zq/a,zq/b,ey}};q,
\frac{qcd}{xab}\right).
\endalign
$$
\endproclaim

\demo{Remark} The conditions for Proposition \thtag{210} are
less severe than in Proposition
\thtag{310}. However, if we want the three ${}_4\vp_3$-series
to be balanced we need the extra conditions
$cd=xab$ and $abq=cde$.
\enddemo

\demo{Proof} As $k\to\infty$, the summand is
${\Cal O}(1)$, so we need $|z|<1$ for absolute
convergence. As $k\to -\infty$ the first
${}_2\vp_1$-series behaves like $C_1 a^{-k}+C_2 b^{-k}$
by \eqtag{330} for $k=1$
and for the second ${}_2\vp_1$-series we use its $q$-Bessel
coefficient behaviour; for $k<0$
$$
\multline
\frac{(q^{k+1};q)_\infty}{(dq^k;q)_\infty}
\, {}_2\vp_1\left( {{dq^k,e}\atop{q^{k+1}}};q,y\right) =
\sum_{j=-k}^\infty \frac{(q^{1+k+j};q)_\infty (e;q)_j}
{(dq^{k+j};q)_\infty (q;q)_j} y^j \\ =
\sum_{p=0}^\infty \frac{(q^{1+p};q)_\infty (e;q)_{p-k}}
{(dq^p;q)_\infty (q;q)_{p-k}} y^{p-k} =
\frac{(q;q)_\infty (e;q)_{-k}}{(d;q)_\infty (q;q)_{-k}}
y^{-k} \, {}_2\vp_1\left( {{d,eq^{-k}}\atop{q^{1-k}}};q,y\right),
\endmultline
$$
which is ${\Cal O}(y^{-k})$ as $k\to-\infty$. So the sum
is absolutely convergent as $k\to-\infty$ for
$|z|>|ya|$ and $|z|>|yb|$.

Let $T$ be the sum, then interchanging summation over $k$
and the summation for the series representation of the $q$-Bessel
coefficient gives
$$
T = \sum_{j=0}^\infty \Bigl( \sum_{k=-j}^\infty
z^k \, {}_2\vp_1\left(
{{a,b}\atop{c}};q,xq^k\right)
\frac{(q^{1+k+j};q)_\infty}{(dq^{k+j};q)_\infty}
\Bigr) \frac{(e;q)_j}{(q;q)_j} y^j .
$$
Consider the inner sum, say $T_j$, for
$j$ fixed. Shifting the summation
parameter and using the series representation
for the ${}_2\vp_1$-series gives
$$
T_j = \sum_{l=0}^\infty z^{-j} \frac{(a,b;q)_l}{(c,q;q)_l}
x^l q^{-jl} \sum_{p=0}^\infty
\frac{(q^{1+p};q)_\infty}{(dq^p;q)_\infty} z^pq^{pl} 
=z^{-j} \frac{(q,dz;q)_\infty}{(d,z;q)_\infty}\, {}_3\vp_2\left(
{{a,b,z}\atop{c,dz}};q,xq^{-j}\right)
$$
using the $q$-binomial theorem. This is only valid
for $|xq^{-j}|<1$, but
since the sum is uniformly convergent on compacta, $T_j$ is
analytic in $x\in\CRp$,
so the result remains valid using the unique analytic
continuation of the ${}_3\vp_2$-series to
$\CRp$.

So this gives
$$
T= \frac{(q,dz;q)_\infty}{(d,z;q)_\infty}
\sum_{j=0}^\infty \frac{(e;q)_j}{(q;q)_j} \Bigl( \frac{y}{z}
\Bigr)^j \, {}_3\vp_2\left(
{{a,b,z}\atop{c,dz}};q,xq^{-j}\right),
\tag\eq{3127}
$$
where we assume that we use the unique analytic continuation
for the ${}_3\vp_2$-series. Since $x\in\CRp$
and $|qcd|<|xab|$ we have for all $j\geq 0$,
$$
\align
&{}_3\vp_2\left(
{{a,b,z}\atop{c,dz}};q,xq^{-j}\right) =
z^j \frac{(a,b,c/z,d,zx,q/zx;q)_\infty}
{(c,dz,a/z,b/z,x,q/x;q)_\infty}
\, {}_3\vp_2\left( {{z, zq/c,q/d}\atop{zq/a,zq/b}};q,
\frac{q^{1+j}cd}{xab}\right)\\
+&a^j\frac{(b,z,c/a,dz/a,ax,q/ax;q)_\infty}
{(c,dz,b/a,z/a,x,q/x;q)_\infty}
\, {}_3\vp_2\left( {{a,aq/c,aq/dz}\atop{aq/b,aq/z}};q,
\frac{q^{1+j}cd}{xab}\right) \\
+&
b^j\frac{(a,z,c/b,dz/b,bx,q/bx;q)_\infty}
{(c,dz,a/b,z/b,x,q/x;q)_\infty}
\, {}_3\vp_2\left( {{b,bq/c,bq/dz}\atop{bq/a,bq/z}};q,
\frac{q^{1+j}cd}{xab}\right)
\endalign
$$
using \eqtag{330} and the theta product identity \eqtag{340}.
Using this in \eqtag{3127} and the $q$-binomial
theorem three times gives the result.
\qed\enddemo

\proclaim{Corollary \thname{332}} Let $S$ be the sum in
Proposition \thtag{310} and assume $ey=q^{1-l}$ for
$l\in\Z$, then the series is absolutely convergent
for $|da|,|db|<|z|<1$. Assuming moreover $abde=cf$ and $fx=dey$
we have that $S$ equals the sum of three balanced
${}_4\vp_3$-series;
$$
\align
S &=
\frac{(q,a,b,c/z,zx,q/zx;q)_\infty}
{(z,c,a/z,b/z,x,q/x;q)_\infty}
\, {}_4\vp_3\left( {{z,zq/c,q^{1-l}/y,d}\atop{zq/a,zq/b,f}};q,
q\right)
\\ &+
\frac{(d,q,b,c/a,zyq^l/a,ax,q/ax,fa/z;q)_\infty}
{(f,yq^l,c,b/a,z/a,x,q/x,da/z;q)_\infty}
\bigl(\frac{z}{a}\bigr)^l
\, {}_4\vp_3\left( {{a,aq/c,aq^{1-l}/yz,ad/z}\atop
{aq/b,aq/z,af/z}};q, q\right)
\\ &+
\frac{(d,q,a,c/b,zyq^l/b,bx,q/bx,fb/z;q)_\infty}
{(f,yq^l,c,a/b,z/b,x,q/x,db/z;q)_\infty}
\bigl(\frac{z}{b}\bigr)^l
\, {}_4\vp_3\left( {{b,bq/c,bq^{1-l}/yz,bd/z}\atop
{bq/a,bq/z,bf/z}};q, q\right).
\endalign
$$
\endproclaim

\demo{Remark} Note that putting $e=q^{1-l}/y$ in the
expression \eqtag{3124} for $S$ gives the same result
using \eqtag{340}. Hence, the expressions for $S$ derived
in Proposition \thtag{310}, \eqtag{3121}, \eqtag{3123} and
\eqtag{3124} remain valid in the case $ey\in q^\Z$.
\enddemo

\demo{Proof} Use \eqtag{3126} in $S$ as in Proposition \thtag{310},
shift the summation parameter $k=n-l$, and apply
Proposition \thtag{331} with
$x\mapsto xq^l$, $e\mapsto f/d$, $d\mapsto yq^l$,
$y\mapsto d$ and use the theta product identity \eqtag{340}
to find an explicit expression for $S$ as a sum of three
${}_4\vp_3$-series for $|da|,|db|<|z|<1$, $|d|<1$,
$x\in\CRp$, $|qcy|<|xab|$. Under the conditions $abde=cf$
and $fx=dey$ we see that these ${}_4\vp_3$-series
are the balanced ${}_4\vp_3$-series as stated.
The condition $|d|<1$ can
be removed by analytic continuation in $d$.
\qed\enddemo

\demo{Proof of Theorem \thtag{210}}
We start off with the general dual transmutation kernel
as in \eqtag{3125} assuming
$|\si|,|\ta|\geq 1$.
This sum is of the type
as considered in Proposition \thtag{310}, and the conditions
in Proposition \thtag{310} lead to $ad=bc$ and $adx/c=by$, or
$ad=bc$ and $x=y$. We put $c=aq^\al$ and $d=bq^\al$, and we
apply Proposition \thtag{310} with $a\mapsto aq^\al/\ta$,
$b\mapsto aq^\al\ta$, $c\mapsto abq^{2\al}$,
$d\mapsto b/\si$, $e\mapsto b\si$, $f\mapsto ab$, $x\mapsto
-bx/a$, $y\mapsto -y$ and $z\mapsto abt$. This gives the
required expression for the dual transmutation kernel under
the conditions $|abq^\al\si\ta|<|abt|<1$ and
$q<ab|t|$, $t\notin a^2q^{2\al+\Zp}$,
$t\notin b^2q^{2\al+\Zp}$.

The previous paragraph deals with the case $\la\in[-1,1]$,
or $\si=e^{i\th}$. For the discrete mass points we first
consider $\si_l=-q^{1-l}/by$, $\l\in\Z$ such that
$|\si_l|>1$. Corollary \thtag{332} then applies,
since $b<1$, $|\si_l|>1$, $q^\al ab<|abt|<1$. From
the remark following Corollary \thtag{332} we see that
the expression for the dual transmutation
kernel remains valid for this set of
discrete mass points of $d\nu(\cdot;a,b;y;q)$. The (possibly
empty) set of discrete mass points of the form $\si_l=aq^l$,
$l\in\Zp$ such that $|\si_l|>1$, make one of the
${}_2\vp_1$-series terminating after using the symmetry
described in Remark \thtag{311}, and the result for
the dual transmutation kernel remains valid. This proves the
first statement of Theorem \thtag{210}.

We now investigate when
$\{t^k \phi_\mu(xq^k;c,d;q)\}_{k\in\Z}\in
{\Cal H}(a,b;y)$. From \eqtag{140} and \eqtag{160}
we see that we need the condition $|t^2ab|<1$ for
convergence as $k\to\infty$. For $k\to -\infty$ we
use the theta product identity \eqtag{340} to find
$$
(ab)^k \frac{(-byq^k/a;q)_\infty}{(-yq^k;q)_\infty}
= a^{2k} \frac{(-q^{1-k}/y;q)_\infty}{(-q^{1-k}a/by;q)_\infty}
\frac{(-by/a,-aq/by;q)_\infty}{(-y,-q/y;q)_\infty}
\tag\eq{3130}
$$
which is ${\Cal O}(a^{2k})$ as $k\to-\infty$.
The asymptotic behaviour of the little $q$-Jacobi
function as $|x|\to\infty$ on a $q$-grid follows from
the expansion \eqtag{180}; for $k\to-\infty$
$$
t^k \phi_\mu(xq^k;c,d;q)=
\cases {\Cal O}\bigl((\frac{t}{c})^k\bigr), &|\ta|=1, \\
{\Cal O}\bigl( (\frac{t}{c\ta})^k\bigr), &|\ta|>1,\ c(\ta;
c,d,x;q)\not=0, \\
{\Cal O}\bigl( (\frac{t\ta}{c})^k\bigr), &|\ta|>1,\ c(\ta;
c,d,x;q)=0.\endcases
\tag\eq{3140}
$$
This implies that for generic $\mu$ we need
$|t^2\ta^{-2}c^{-2}a^2|>1$ for convergence of
\eqtag{160} as $k\to-\infty$. We conclude
$$
\bigl\vert \frac{c\ta}{a}\bigr| < |t| <\frac{1}{\sqrt{|ab|}}
\Longrightarrow
\{t^k \phi_\mu(xq^k;c,d;q)\}_{k\in\Z}\in
{\Cal H}(a,b;y).
\tag\eq{3150}
$$
Under the assumption of \eqtag{3150} the $L^2$-theory for
the little $q$-Jacobi function transform implies
that the sum in \eqtag{3125} converges in
$L^2\bigl(d\nu(\cdot;a,b;y;q)\bigr)$. In the special
case $c=aq^\al$, $d=bq^\al$, $x=-bx/a$, the conditions
\eqtag{3150} are implied by $abq^\al|\ta|<ab|t|<1$, since
$0<ab<1$. Hence, $P_t(\cdot,\mu;q^\al;a,b)
\in L^2(d\nu(\cdot;a,b;y;q))$ for $q^\al|\ta|<|t|<1/\sqrt{ab}$.

Assuming $|q^\al\ta|<|t|<1/\sqrt{ab}$ we find from \eqtag{170}
$$
t^k\,\phi_\mu(yq^k;aq^\al,bq^\al;q) = \int_\R
\phi_\la(yq^k;a,b;q) \, P_t(\la,\mu;q^\al;a,b)
\, d\nu(\la;a,b;y;q).
$$
Taking linear combinations shows that for $v=\{v_k\}_{k\in\Z}$
with only finitely many non-zero coefficients, we have
$$
{\Cal F}_{aq^\al,bq^\al,y}\bigl[\de_t \, v](\mu)
= \int_\R {\Cal F}_{a,b,y}\bigl[k\mapsto q^{2\al k}v_k\bigr](\la)
 \, P_t(\la,\mu;q^\al;a,b)
\, d\nu(\la;a,b;y;q),
\tag\eq{3160}
$$
using the notation as in \eqtag{181}.
Again using the $L^2$-theory of \eqtag{170}, \eqtag{3160}
remains
valid for $\{q^{2k\al}v_k\}_k\in {\Cal H}(a,b;y)$, since
this makes the integrand  integrable with respect to
$d\nu(\cdot;a,b;y;q)$. Now we take $v_k=s^kt^{-k}
\phi_\nu(yq^k;aq^{\al+\be},bq^{\al+\be};q)$, which
gives the product formula. This is valid for
$|q^{\al+\be}\rho|<|q^{2\al}s/t|<1/\sqrt{ab}$ by
\eqtag{3150}, where $\nu=\hf(\rho+\rho^{-1})$ with $|\rho|\geq 1$.
Note that these two conditions on
$s$ and $t$ imply $|abq^{2\al+\be}\rho\ta|<|abq^{2\al}s|<1$,
which are precisely the conditions for the absolute convergence
of the left hand side of \eqtag{3160}, cf. \eqtag{3125}.
\qed\enddemo

\head \newsec The transmutation kernel \endhead

In this section we prove Theorem \thtag{230}. The results in
this section give another point of view to Gasper's results
\cite{\GaspCM} on $q$-analogues of Erd\'elyi's fractional integrals,
see also \S 5.

Rewriting Gasper's $q$-analogue \cite{\GaspCM, (1.8)} of
Erd\'elyi's fractional integral gives
$$
\multline
{}_2\vp_1\left({{ar\si,ar/\si}\atop{abrs}};q,
-\frac{byq^l}{ar}\right) =
\frac{(ab,rs;q)_\infty}{(q,abrs;q)_\infty}\sum_{k=0}^\infty
(ab)^k \frac{(q^{k+1},-byq^{k+l}/a;q)_\infty}
{(rsq^k,-byq^{k+l}/ar;q)_\infty} \\ \times
\, {}_3\vp_2\left( {{q^{-k},r,ar/b}\atop{rs, -arq^{1-l-k}/by}};q,q
\right) \, {}_2\vp_1\left({{a\si,a/\si}\atop{ab}};q,
-\frac{byq^{l+k}}{a}\right),
\endmultline
\tag\eq{610}
$$
for $|rs|<1$, $|ab|<1$.
In \S 5, see Theorem \thtag{240}(ii),
we give an alternative derivation of \eqtag{610} using
intertwining operators for the second order $q$-difference
operator $L$ as in \eqtag{150}. Indeed, as remarked
in \S 2, \eqtag{250}
is equivalent to \eqtag{610}.
We can also prove \eqtag{610} using
the connection coefficient formula
$$
(b\si,b/\si;q)_m = (ab,\frac{b}{a};q)_m \sum_{k=0}^m
\frac{(q^{-m};q)_k q^k}{(q,ab,q^{1-m}a/b;q)_k}
(a\si,a/\si;q)_k,
\tag\eq{612}
$$
which is the $q$-Saalsch\"utz
summation formula \cite{\GaspR, (II.12)}.
So \eqtag{610} links two little $q$-Jacobi functions with different
parameters sets, but with related argument.
Without this condition on the arguments we have not been able to
find a simple expression for the kernel, even if we
allow the summation
parameter to run over $\Z$ instead of $\Zp$.

Note that \eqtag{610} may be viewed as a connection
coefficient problem for two sets of orthogonal functions. Having in
mind the limit transition of Askey-Wilson polynomials to little
$q$-Jacobi functions, see \cite{\KoelSNATO,
\S\S 2.3, 4.1, 6.1}, we
see that \eqtag{610} can be viewed as the limit case of the
connection coefficient problem for Askey-Wilson polynomials, see
\cite{\AskeW, \S 6}.

In case $s=1$ the ${}_3\vp_2$-series in \eqtag{610}
reduces to a terminating ${}_2\vp_1$-series that can be summed by
the $q$-Chu-Vandermonde summation
\cite{\GaspR, (II.6)} yielding
$$
\multline
{}_2\vp_1\left({{ar\si,ar/\si}\atop{abr}};q,
-\frac{byq^l}{ar}\right) =
\frac{(ab,r,-yq^l;q)_\infty}{(q,abr,-byq^l/ar;q)_\infty}
\\ \sum_{k=0}^\infty
(ab)^k \frac{(q^{k+1},-byq^{k+l}/a;q)_\infty}
{(rq^k,-yq^{k+l};q)_\infty}
\, {}_2\vp_1\left({{a\si,a/\si}\atop{ab}};q,
-\frac{byq^{l+k}}{a}\right),
\endmultline
\tag\eq{620}
$$
which is equivalent to, by \cite{\GaspR, (1.4.6)},
$$
{}_2\vp_1\left({{b\si,b/\si}\atop{abr}};q,
-yq^l\right) =
\frac{(ab,r;q)_\infty}{(q,abr;q)_\infty}
\sum_{k=0}^\infty
(ab)^k \frac{(q^{k+1};q)_\infty}
{(rq^k;q)_\infty}
\, {}_2\vp_1\left({{b\si,b/\si}\atop{ab}};q,
-yq^{l+k}\right),
\tag\eq{630}
$$
valid for $y\in\CRp$. This can be proved easily using
the $q$-binomial theorem. See also
\thetag{5.7} 
and \thetag{5.13} 
for similar results.
Another interesting case of \eqtag{610} is $r=1$, which
reduces the ${}_3\vp_2$-series to $1$. This gives \eqtag{630}
after renaming. In \S 5, see Theorem \thtag{240},
we show that \eqtag{610} can be derived
from \eqtag{630}.

\demo{Proof of Theorem \thtag{230}}
Fix $y$ and define $u=\{u_k\}_{k\in\Z}$ by
$$
u_k = (ab)^{-l}\frac{(ab,rs;q)_\infty}{(q,abrs;q)_\infty}
\frac{(q^{k-l+1},-yq^k;q)_\infty}
{(rsq^{k-l},-byq^k/ar;q)_\infty}
\, {}_3\vp_2\left( {{q^{l-k},r,ar/b}\atop{rs, -arq^{1-k}/by}};q,q
\right)
$$
with the convention $u_k=0$ for $k<l$. Using
\cite{\GaspR, (III.12)} the ${}_3\vp_2$-series can be written
as
$$
{}_3\vp_2\left( {{q^{l-k},r,ar/b}\atop{rs, -arq^{1-k}/by}};q,q
\right) =  \frac{(-yq^l;q)_{k-l}}{(-byq^l/ar;q)_{k-l}}
\, {}_3\vp_2\left( {{q^{l-k},ar/b,s}\atop{rs, -yq^l}};q,
-\frac{byq^k}{a}\right),
$$
so that $u_k={\Cal O}(1)$ as $k\to\infty$.
It follows that $u\in {\Cal H}(a,b;y)$.
Then \eqtag{610} states
that the little $q$-Jacobi transform of $u$, see \eqtag{170},
is
$$
\bigl( {\Cal F}_{a,b,y}u\bigr) (\la)=
\phi_\la( yq^l/s;ar,bs;q) =
{}_2\vp_1\left({{ar\si,ar/\si}\atop{abrs}};q,
-\frac{byq^l}{ar}\right), \qquad \la=\hf(\si+\si^{-1}).
$$
By the inversion formula of \eqtag{170} we find the result
for the transmutation kernel as in Theorem \thtag{230}.

Note that we can write
$$
\sum_{l=-\infty}^\infty u_l (abrs)^l \frac{(-byq^l/ar;q)_\infty}
{(-yq^l/s;q)_\infty} \, P_{k,l}(a,b,y;r,s)
= \bigl( {\Cal F}_{a,b,y}^{-1} {\Cal F}_{ar,bs,y/s} u\bigl)_k
\tag\eq{660}
$$
valid for $u=\{ u_k\}_k$ having only finitely many non-zero
coefficients.
In particular, choosing $u_l=P_{l,p}(ar,bs,y/s;t,u)$ for
$l\leq k$ and $u_l=0$ for $l>k$ under
the assumptions $u,t>0$, $|ut|<1$ we find from
\eqtag{660}
$$
\multline
\sum_{l=-\infty}^\infty P_{k,l}(a,b,y;r,s)
P_{l,p}(ar,bs,y/s;t,u)\, (abrs)^l \frac{(-byq^l/ar;q)_\infty}
{(-yq^l/s;q)_\infty}
\\ = \bigl( {\Cal F}_{a,b,y}^{-1} {\Cal F}_{ar,bs,y/s}
\bigl[ l\mapsto \bigl(
{\Cal F}_{ar,bs,y/s}^{-1}[\la\mapsto
\phi_\la(\frac{yq^p}{su};art,bsu;q)]_l\bigr] \bigl)_k \\ =
\bigl( {\Cal F}_{a,b,y}^{-1}[\la\mapsto
\phi_\la(\frac{yq^p}{su};art,bsu;q)]\bigl)_k =
P_{k,p}(a,b,y;rt,su),
\endmultline
\tag\eq{670}
$$
which is the product formula.
\qed\enddemo

If we plug the explicit expression for the transmutation
kernel into the product formula, we obtain, after
relabelling, the following expression. For $k\in\Z$, $p\in\Zp$,
$r,s,t,u>0$ with $rs<1$, $tu<1$ we have
$$
\multline
\sum_{l=0}^p
\frac{(tu;q)_{p-l}}{(q;q)_{p-l}} \, {}_3\vp_2\left(
{{q^{l-p},t,art/bs}\atop{tu, -artq^{1+l-k}/by}};q,q\right)
\frac{(rs;q)_l}{(q;q)_l} \, {}_3\vp_2\left(
{{q^{-l},r,ar/b}\atop{rs, -arq^{1-k}/by}};q,q\right)
\\ \times (rs)^{p-l}t^l
\frac{(-arq^{1-k}/by;q)_l}{(-artq^{1-k}/by;q)_l}
=  \frac{(rstu;q)_p}{(q;q)_p} \, {}_3\vp_2\left(
{{q^{-p},rt,art/b}\atop{rstu, -artq^{1-k}/by}};q,q\right).
\endmultline
\tag\eq{685}
$$
This product formula is equivalent to Gasper \cite{\GaspCM, (1.7)},
and, as remarked in \cite{\GaspCM}, \eqtag{685} implies \eqtag{610}
in the limit $p\to\infty$. Here we have shown that
the converse is also valid, \eqtag{610} implies \eqtag{685}
using the little $q$-Jacobi function transform.

\head \newsec Intertwining properties
\endhead

In this section we prove Theorem \thtag{240}, and as
a motivation we start by giving a Darboux factorisation of
the second order $q$-difference operator
$L^{(a,b)}$ or ${\Cal L}^{(a,b)}$.

The backward $q$-derivative
operator is $B_q=M_{1/x}(1-T_q^{-1})$, where $M_g$ is the
operator of multiplication by $g$;
$\bigl( M_gf\bigr)(x)=g(x)f(x)$, and $T_qf(x)=f(qx)$,
as introduced in \S 1. It is straightforward
to check that
$$
\bigl( B_q\phi_\la(\cdot;a,b;q)\bigr)(x)
= \frac{b(1-a\si)(1-a/\si)}{qa(1-ab)}
\, \phi_\la(x;aq,b;q).
\tag\eq{501}
$$
Considering ${\Cal H}(a,b;y)$ as an $L^2$-space with discrete
weights $(ab)^k (-byq^k/a;q)_\infty/(-yq^k;q)_\infty$ at
the point $yq^k$, $k\in\Z$, we look at $B_q$ as a
(densely defined unbounded) operator
from ${\Cal H}(a,b;y)$ to ${\Cal H}(aq,b;y)$. Its adjoint,
up to a constant depending only on $y$,
is given by
$$
A(a,b) = M_{1+bx/aq} - ab M_{1+x}T_q,
\tag\eq{502}
$$
and it's a straightforward calculation to show that
$$
\bigl( A(a,b) \phi_\la(\cdot;aq,b;q)\bigr)(x) =
(1-ab)\, \phi_\la(x;a,b;q)
\tag\eq{503}
$$
and that
$-b L^{(a,b)} = aq A(a,b) \circ B_q$, with the
notation as in \eqtag{150}.
Since
$B_q$ and $A(a,b)$ are triangular with respect
to the standard orthogonal basis of Dirac delta's at $yq^k$
of ${\Cal H}(a,b;y)$, this means that we have
a Darboux factorisation of $L^{(a,b)}$. Also,
$-b(L^{(aq,b)}+(1-q)(1-qa^2))=aq^2 B_q\circ A(a,b)$, from
which we deduce $B_q\circ L^{(a,b)}=L^{(aq,b)}\circ B_q$
and $L^{(a,b)}\circ A(a,b)=A(a,b)\circ L^{(aq,b)}$.
It is the purpose of this section to generalise these
intertwining properties to arbitrary complex powers of $B_q$.

Introduce the operator $W_\nu$, $\nu\in\C$, acting on
functions defined on $[0,\infty)$ by
$$
\bigl( W_\nu f\bigr)(x) = x^\nu \sum_{l=0}^\infty
f(xq^{-l}) q^{-l\nu} \frac{(q^\nu;q)_l}
{(q;q)_l}, \qquad x\in[0,\infty),
\tag\eq{510}
$$
assuming that the
infinite sum is absolutely convergent if $\nu\notin -\Zp$.
So we want
$f$ sufficiently decreasing on a $q$-grid tending to
infinity, e.g. $f(xq^{-l}) ={\Cal O}(q^{l(\nu+\ep)})$
for some $\ep>0$. Note that for $\nu\in\Z_{\leq 0}$ the sum
in \eqtag{510} is finite and $W_0=\text{Id}$ and $W_{-1}=B_q$.

This operator is a $q$-analogue of the Weyl fractional
integral operator as used in \cite{\KoorAM, \S 3},
\cite{\KoorJF, \S 5.3} for the Abel transform.
With the notation
$$
\int_a^\infty f(t)\, d_qt = a\sum_{k=0}^\infty f(xq^{-k})q^{-k}
$$
for the $q$-integral, cf. \eqtag{165}, we see that for
$n\in\N$ the operator $W_n$ is an iterated $q$-integral;
$$
\bigl( W_n f\bigr)(x) = \int_x^\infty \int_{x_1}^\infty \ldots
\int_{x_{n-1}}^\infty f(x_n)\, d_qx_nd_qx_{n-1}\ldots d_qx_1.
\tag\eq{512}
$$

In the following lemma we collect some results on
$W_\nu$, where we use the function space
$$
{\Cal F}_\rho = \{ f\colon [0,\infty)\to \C \mid\,
|f(xq^{-l})| = {\Cal O}(q^{l\rho})
,\ l\to\infty,\ \forall x\in (q,1]\},\qquad \rho>0.
\tag\eq{515}
$$
Recall that ${\Cal L}^{(a,b)}$ is defined in
\eqtag{155}.

\proclaim{Lemma \thname{520}} Let
$\nu,\mu\in\C\backslash\Z_{\leq 0}$.
\item{\rm (i)} $W_\nu$ preserves the
space of compactly supported functions,
\item{\rm (ii)} $W_\nu\colon
{\Cal F}_\rho\to {\Cal F}_{\rho-\Re\nu}$ for $\rho>\Re\nu>0$,
\item{\rm (iii)} $W_\nu\circ W_\mu = W_{\nu+\mu}$
on ${\Cal F}_\rho$ for $\rho>\Re(\mu+\nu)>0$,
\item{\rm (iv)} $W_\nu\circ B_q = B_q\circ W_\nu = W_{\nu-1}$
on ${\Cal F}_\rho$ for $\rho>\Re\nu-1>0$,
and $B_q^n\circ W_n = \id$ for $n\in\N$ on ${\Cal F}_\rho$
for $\rho>n$,
\item{\rm (v)} ${\Cal L}^{(aq^{-\nu},b)}\circ W_\nu =
W_\nu \circ {\Cal L}^{(a,b)}$, valid for compactly
supported functions.
\endproclaim

\demo{Remark} It follows from (iii) that
$W_{-n}=B_q^n$, $n\in\N$, and $W_0=\id$.
\enddemo

\demo{Proof} The first statement is immediate from \eqtag{510}.
For (ii) we use that for $f\in{\Cal F}_\rho$ and $x\in(q,1]$
we have
$$
| W_\nu f(xq^{-k})| \leq M\sum_{l=0}^\infty q^{(k+l)\rho}
q^{-(k+l)\Re\nu} \frac{(q^{\Re\nu};q)_l}{(q;q)_l} =
M q^{k(\rho-\Re\nu)}
\frac{(q^\rho;q)_\infty}{(q^{\rho-\Re\nu};q)_\infty}
$$
by the $q$-binomial theorem for
$\rho>\Re\nu$.
The third statement is a consequence of interchanging
summations, valid for $f\in{\Cal F}_\rho$, $\rho>\Re(\mu+\nu)$,
and
$$
\sum_{k+l=p} \frac{(q^\mu;q)_k (q^\nu;q)_l}
{(q;q)_k(q;q)_l} q^{-(l+k)\mu-l\nu} =
q^{-p(\mu+\nu)} \frac{(q^{\mu+\nu};q)_p}
{(q;q)_p},
$$
which is the $q$-Chu-Vandermonde summation formula
\cite{\GaspR, (1.5.2)}. For (iv) we note that
$B_q\colon {\Cal F}_\rho\to {\Cal F}_{\rho+1}$, then
the first statement of (iv) is a simple calculation
involving $q$-shifted factorials, which reduces
the second statement of (iv) to verifying the easy case $n=1$.
For (v) recall \eqtag{155},
so that ${\Cal L}^{(aq^{-\nu},b)}(W_\nu f)(x)$ and  $W_\nu ({\Cal
L}^{(a,b)}f)(x)$ involve the values $f(xq^{-k})$,
$k+1\in \Zp$. A straightforward calculation
using $q$-shifted factorials shows that the coefficients
of $f(xq^{-k})$ in
${\Cal L}^{(aq^{-\nu},b)}(W_\nu f)(x)$ and  $W_\nu ({\Cal
L}^{(a,b)}f)(x)$ are equal.
\qed\enddemo

The asymptotically free solution $\Phi_\si(yq^k;a,b;q)
\in{\Cal F}_\rho$ for $q^\rho>|a\si|$ as follows from
\eqtag{180}. A calculation using the $q$-binomial formula gives,
cf. \eqtag{630},
$$
\bigl(W_\nu \Phi_\si(\cdot;a,b;q)\bigr)(yq^k)=
y^\nu \frac{(a\si;q)_\infty}
{(aq^{-\nu}\si;q)_\infty}
\Phi_\si(yq^k;aq^{-\nu},b;q),
\tag\eq{540}
$$
for $|a\si|<q^\nu$ in accordance with Lemma \thtag{520}(v).
Note that \eqtag{540} is a $q$-analogue of Bateman's formula,
cf. \cite{\GaspCM}, \cite{\KoorAM}.

\proclaim{Lemma \thname{542}} Define the operator
$$
S(a,b) = M_{\frac{(-x;q)_\infty}{(-bx/a;q)_\infty}}
\circ T_{b/a},\qquad T_{b/a}f(x) = f(\frac{b}{a}x),
$$
then $S(a,b)^{-1}\circ {\Cal L}^{(a,b)} \circ
S(a,b) = {\Cal L}^{(b,a)}$.
In particular, $\tilde W_\nu^{(a,b)} = S(a,bq^{-\nu})\circ
W_\nu\circ S(a,b)^{-1}$ satisfies the intertwining property
${\Cal L}^{(a,bq^{-\nu})}\circ \tilde W_\nu^{(a,b)} =
\tilde W_\nu^{(a,b)} \circ {\Cal L}^{(a,b)}$.
\endproclaim

Note that $S(a,b)^{-1}=S(b,a)$ and that $S(a,b)\colon
{\Cal H}(b,a;yb/a)\to {\Cal H}(a,b;y)$ is an isometric
isomorphism.
For $f\in{\Cal F}_\rho$ we see that
$\bigl(S(a,b)f\bigr)(xq^{-l}) = {\Cal O}( |a/b|^lq^{l\rho})$
from \eqtag{340}, so that
$S(a,b)f\in {\Cal F}_{\rho + \ln(|a/b|)/\ln q}$.

\demo{Proof} It follows from \eqtag{155} that
$$
\multline
{\Cal L}^{(a,b)}\bigl(x\mapsto \frac{(-x;q)_\infty}
{(-bx/a;q)_\infty}f(x)\bigr)(x) = \\ \frac{(-x;q)_\infty}
{(-bx/a;q)_\infty} \Bigl( \frac{b}{2}(1+\frac{a}{bx})f(qx)
+ \frac{1}{2b}(1+\frac{q}{x})f(xq^{-1}) - \hf
(\frac{a}{x}+\frac{q}{bx})f(x)\Bigr)
\endmultline
$$
and the term in parentheses can be written as
$T_{b/a}\circ {\Cal L}^{(b,a)}\circ T_{a/b}$ applied to $f$.
The second statement then follows from Lemma \thtag{520}(v).
\qed\enddemo

It follows directly from \eqtag{140}, \eqtag{180} and
\cite{\GaspR, (1.4.6)} that
$$
\aligned
\bigl( S(a,b)\phi_\la(\cdot;b,a;q)\bigr)(x)&= \phi_\la(x;a,b;q), \\
\bigl( S(a,b)\Phi_\si(\cdot;b,a;q)\bigr)(x)&= \Phi_\si(x;a,b;q).
\endaligned
\tag\eq{544}
$$

\demo{Proof of the first statement of Theorem \thtag{240}}
It follows from Lemma \thtag{520}(v) and Lemma \thtag{542}
that the operator
$$
\align
W_{\nu,\mu}(a,b) &= \tilde W_\mu^{(aq^{-\nu},b)}\circ W_\nu \\
&= S(aq^{-\nu},bq^{-\mu})\circ W_\mu\circ
S(b,aq^{-\nu})\circ W_\nu
\endalign
$$
satisfies the required interwining property. For $f\in {\Cal
F}_\rho$ with $\rho>\Re\nu$  we can interchange summations,
which leads
to the sum with a terminating ${}_3\vp_2$ as kernel. Note that
the ${}_3\vp_2$-series in the kernel of $W_{\nu,\mu}(a,b)$
behaves as
$$
{}_2\vp_1\left( {{q^{-\mu},-q^{1+\mu-\nu}a/bx}\atop{-q^{1+\mu}/x}};
q, q^{\nu-\mu}\frac{b}{a}\right)
$$
as $p\to\infty$.

The statement for the action on $\Phi_\si(\cdot;a,b;q)$ follows
immediately from \eqtag{540} and \eqtag{544}.
\qed\enddemo

In order to prove the remaining half of Theorem \thtag{240} we
take appropriate adjoints of the previous construction.
Consider $W_\nu$, $\nu\in\C\backslash\Z_{\leq 0}$, as a
densely defined unbounded operator from
${\Cal H}(aq^\nu,b;y)$ to ${\Cal H}(a,b;y)$ and define
$R_\nu^{(a,b)}$ as its adjoint, so
$$
\langle R^{(a,b)}_\nu f, g\rangle_{{\Cal H}(aq^\nu,b;y)} =
\langle f, W_\nu g\rangle_{{\Cal H}(a,b;y)}
\tag\eq{545}
$$
for all compactly supported functions $g$, cf. Lemma \thtag{520}(i).
Here we use the identification of ${\Cal H}(a,b;y)$ as a weighted
$L^2$-space on a discrete set, see \S 1. A $q$-integration by parts
shows
$$
\bigl(R_\nu^{(a,b)} f\bigr)(yq^p)= y^\nu
\frac{(-byq^p/a;q)_\infty}{(-byq^{p-\nu}/a;q)_\infty}
\sum_{l=0}^\infty f(yq^{p+l}) (ab)^l
\frac{(q^\nu,-yq^p;q)_l}{(q,-byq^p/a;q)_l}.
\tag\eq{560}
$$
Now define, for functions $f$, the operator
$$
\bigl(A_\nu^{(a,b)} f\bigr)(x)=
\frac{(-bx/a;q)_\infty}{(-bxq^{-\nu}/a;q)_\infty}
\sum_{l=0}^\infty f(xq^l) (ab)^l
\frac{(q^\nu,-x;q)_l}{(q,-bx/a;q)_l},
\tag\eq{570}
$$
so that $A_\nu^{(a,b)}\Big\vert_{{\Cal
H}(a,b;y)}=y^{-\nu}R_\nu^{(a,b)}$. Note that
$A_\nu^{(a,b)}$ is well-defined for bounded functions
assuming $|ab|<1$.
Recall that the dense
domain of finite linear combinations of the basis
vectors for ${\Cal L}^{(a,b)}$ corresponds to the
functions compactly supported in $(0,\infty)$.

\proclaim{Lemma \thname{550}}
${\Cal L}^{(aq^{\nu},b)}\circ A_\nu^{(a,b)} =
A_\nu^{(a,b)} \circ {\Cal L}^{(a,b)}$
on the space of functions
compactly supported in $(0,\infty)$. Moreover,
$$
\bigl(A_\nu^{(a,b)}\phi_\la(\cdot;a,b;q)\bigr)(x) =
 \frac{(abq^\nu;q)_\infty}{(ab;q)_\infty}
\phi_\la(x;aq^\nu,b;q).
$$
Defining $\tilde A_\nu^{(a,b)} = S(a,bq^\nu)\circ
A_\nu^{(b,a)}\circ S(b,a)$ we have
${\Cal L}^{(a,bq^\nu)}\circ \tilde A_\nu^{(a,b)} =
\tilde A_\nu^{(a,b)} \circ {\Cal L}^{(a,b)}$, and
$$
\bigl(\tilde A_\nu^{(a,b)}\phi_\la(\cdot;a,b;q)\bigr)(x) =
\frac{(abq^\nu;q)_\infty}{(ab;q)_\infty}
\phi_\la(x;a,bq^\nu;q).
$$
\endproclaim

\demo{Proof} Note that \eqtag{560} and \eqtag{570} show that
the operators $R_\nu^{(a,b)}$ and $A_\nu^{(a,b)}$ preserve
the space of functions compactly supported in $(0,\infty)$.
The intertwining property for $R_\nu^{(a,b)}$
follows from \eqtag{545} and
Lemma \thtag{520}, and hence for $A_\nu^{(a,b)}$.

To calculate the action of $A_\nu^{(a,b)}$ on the
little $q$-Jacobi function we use
\cite{\GaspR, (1.4.6)} to write
$$
\phi_\la(x;a,b;q) = \frac{(-x;q)_\infty}{(-bx/a;q)_\infty}
\, {}_2\vp_1\left( {{b\si,b/\si}\atop{ab}};q,-x\right).
\tag\eq{575}
$$
Using this in \eqtag{570}, interchanging summations,
which is easily justified for $|x|<1$,
and using the $q$-binomial theorem gives
$$
\bigl(A_\nu^{(a,b)}\phi_\la(\cdot;a,b;q)\bigr)(x) =
\frac{(abq^\nu,-x;q)_\infty}
{(ab,-bxq^{-\nu}/a;q)_\infty}
\, {}_2\vp_1\left( {{b\si,b/\si}\atop{abq^\nu}};q, -x\right)
\tag\eq{580}
$$
and using \eqtag{575} again gives the result for $|x|<1$.
The general case follows by analytic continuation in $x$,
see \eqtag{180}, \eqtag{330},
since the convergence in \eqtag{570} for $f$ the little $q$-Jacobi
function is uniform on compact sets for $x$.

The
statements for $\tilde A_\nu^{(a,b)}$ follow from the
corresponding
statements for $A_\nu^{(a,b)}$ and Lemma \thtag{542} and
\eqtag{544}.
\qed\enddemo

\demo{Proof of the second statement of Theorem \thtag{240}}
Define
$$
\align
A_{\nu,\mu}(a,b) &= \tilde A^{(aq^\nu,b)}_\mu\circ A_\nu^{(a,b)}\\
&= S(aq^\nu,bq^\mu)\circ A_\mu^{(b,aq^\nu)}\circ
S(b,aq^\nu)\circ A_\nu^{(a,b)}
\endalign
$$
then it follows from Lemma \thtag{550} that the intertwining
property is valid. The action on a function $f$ can be calculated
and for $f$ compactly suppported in $(0,\infty)$
we find the explicit result
with the ${}_3\vp_2$-series as kernel.
We can extend the result to bounded $f$ if we require $\nu>0$.

The action of $A_{\nu,\mu}(a,b)$ on the little $q$-Jacobi
function follows from Lemma \thtag{550}.
\qed\enddemo

\head \newsec Quantum group theoretic interpretation
\endhead

The quantised universal enveloping algebra $U_q(\frak{su}(1,1))$
has representations in $\Hp$ for
the discrete series representations,
and in $\H$ for the principal unitary series, the complementary
series and strange series representations. For the harmonic
analysis the so-called twisted primitive elements,
as analogues of self-adjoint Lie algebra elements, play an
important role, and in each of these representations they
give rise to a three-term recurrence relations in which
there is essentially one degree of freedom. In the
discrete series representations the three-term recurrrence
relations can be solved in terms of Al-Salam and Chihara
polynomials \cite{\KoelVdJSIAM}, \cite{\RoseCM},
and in the other series in terms of little $q$-Jacobi functions
\cite{\KoelSsu, \S 6}.

The transition of the generalised
basis of eigenvectors of two different twisted primitive
elements in a strange series representation of
$U_q(\frak{su}(1,1))$
is given by the dual
transmutation kernel of Theorem \thtag{210}.
For the complementary series and principal unitary series
we can deduce the corresponding dual transmutation kernel
from Proposition \thtag{310} in a similar way
using other specialisations, cf. \cite{\KoelSsu, \S 6, App.~A}.
For the discrete series representation
we refer to \cite{\RoseCM} and \cite{\KoelVdJCA}.


\enddocument